\documentclass[10pt]{amsart}
\usepackage{latexsym}
\usepackage[dvips]{color}
\usepackage{float}
\usepackage{afterpage} 
\usepackage [pdftex]{graphicx}
\usepackage{epsfig}
\usepackage{amsmath}
\usepackage{amsthm}
\usepackage{amsfonts}
\usepackage{amssymb}
\usepackage{setspace}
\usepackage{subfigure}
\usepackage{tikz}
\usetikzlibrary{arrows}
\usetikzlibrary{arrows.meta}
\usetikzlibrary{matrix}
\usepackage{subfigure}

\newtheorem{theorem}{Theorem}
\newtheorem{definition}{Definition}
\newtheorem{proposition}{Proposition}

\newtheorem{lemma}{Lemma}

\newcommand{\Aff}{\operatorname{Aff}_{+}(\mathbb R)}

\newcommand{\VertexGap}{20mm}
\newcommand{\LaneGap}{8mm}
\newcommand{\VertexSize}{4mm}

\tikzset{
	vertex/.style={
		draw=black!75,
		circle,
		fill=white,
		inner sep=2pt,
		minimum size=\VertexSize,
		line width=.6pt,
		font=\small
	},
	forward/.style={
		draw=white!50!black,
		line width=.45pt,
		-{Stealth[length=2.2mm,width=1.5mm]}
	},
	backward/.style={
		draw=white!50!black,
		line width=.45pt,
		-{Stealth[length=2.2mm,width=1.5mm]}
	},
	edge label/.style={
		fill=white,
		inner xsep=3pt,
		inner ysep=2pt,
		font=\small
	}
}

\makeatletter

\newcommand{\TransformationDiagram}[2]{%
	\begin{tikzpicture}
		
		\def\nodelabels{#1}
		\def\edges{#2}
		\pgfmathtruncatemacro{\N}{0}
		\foreach \label in \nodelabels {
			\pgfmathtruncatemacro{\N}{\N+1}
		}
		
		\pgfmathsetmacro{\pitch}{\VertexGap/1cm}
		\pgfmathsetmacro{\lane}{\LaneGap/1cm}
		\pgfmathsetmacro{\radius}{\VertexSize/2/1cm}
		
		\pgfmathsetmacro{\firstlane}{\radius+0.8}
		
		\foreach \label [count=\i from 0] in \nodelabels {
			\node[vertex] (v\i) at (\i*\pitch,0) {\label};
		}
		
		\foreach \i/\j/\label in \edges {
			
			\ifnum\i<\j
			
			\pgfmathsetmacro{\height}
			{\firstlane+(\j-\i-1)*\lane}
			\pgfmathsetmacro{\control}{(4*\height-\radius)/3}
			
			\draw[forward]
			(v\i.north)
			.. controls
			(\i*\pitch+0.38*\pitch,\control)
			and
			(\j*\pitch-0.38*\pitch,\control)
			..
			(v\j.north);
			
			\node[edge label,text=white!30!black]
			at ({(\i+\j)*\pitch/2},\height)
			{\label};
			
			\else
			\ifnum\i>\j
			
			\pgfmathsetmacro{\height}
			{\firstlane+(\i-\j-1)*\lane}
			\pgfmathsetmacro{\control}{(4*\height-\radius)/3}
			
			\draw[backward]
			(v\i.south)
			.. controls
			(\i*\pitch-0.38*\pitch,-\control)
			and
			(\j*\pitch+0.38*\pitch,-\control)
			..
			(v\j.south);
			
			\node[edge label,text=white!30!black]
			at ({(\i+\j)*\pitch/2},-\height)
			{\label};
			
			\fi
			\fi
		}
		
	\end{tikzpicture}%
}

\begin{document}

\title[Non-abelian homometry on the positive affine line group]{On simultaneous left- and right-homometric sets on the positive affine line group.}
\author{Alexandre Popoff}
\address{Independent Researcher}
\email{al.popoff@free.fr}


\begin{abstract}
On a non-abelian group $G$, the homometry problem comes in two flavors, depending on whether differences are computed on the left or on the right. A natural question is to determine whether simultaneously left- and right-homometric subsets of $G$ exist. This question has been answered positively in the case of dihedral groups by Genuys and Popoff. We here prove that no such subsets exist on the positive affine line group, if we impose that the subgroups generated by their left- and right-differences be non-abelian.
\end{abstract}

\maketitle
\section{Introduction}

This contribution deals with the theory of homometric sets, which appear in various forms in multiple domains: phase retrieval in crystallography \cite{Rosenblatt1984,GrimmBaake2008}, in optimization problems such as the `turnpike problem' on the integer line or the `beltway problem' on a finite cyclic group \cite{Lemke2003}, as the so-called `Z-relation' in mathematical musicology \cite{Forte1973,Jedrzejewski2013}. Given an abelian group $G$, let $X$ be a subset of $G$ and call $D_X = \{gh^{-1}, g \in X, h\in X\}$ the multiset of differences (or distances) between elements of $X$.

\begin{definition}
	Two subsets $A$ and $B$ of $G$ are said to be \textit{homometric} if their multisets of differences $D_A$ and $D_B$ are identical.
\end{definition}

Evidently, if $A$ and $B$ are trivially homometric if $B$ is related to $A$ by a translation, i.e. $B=gA$ for $g \in G$, or if $B$ is obtained from $A$ by the map $x \to x^{-1}$. It is therefore understood, when speaking about homometric sets, that only \textit{non-trivial} homometry is retained. On $G=\mathbb{Z}$, the smallest pair of representatives is known to be $A=\{0,1,2,6,8,11\}$ and $B=\{0,1,6,7,9,11\}$, giving the distance multiset $\{0,0,0,0,0,0,1,1,2,2,3,4,5,5,6,6,7,8,9,10,11\}$. On $G=\mathbb{Z}_{12}$, a very well-known example from Forte \cite{Forte1973} is the pair $A=\{0,1,4,6\}$ and $B=\{0,1,3,7\}$, corresponding to the pitch-class sets $\{C,C\sharp,E,F\sharp\}$ and $\{C,C\sharp,E\flat,G\}$ in musical notation (assuming $C=0$).

If $G$ is non-abelian, homometry now comes in two flavors. More precisely, for $X$ a subset of $G$ call $L_X = \{gh^{-1}, g \in X, h\in X\}$ and $R_X = \{g^{-1}h, g \in X, h\in X\}$.

\begin{definition}
	Two subsets $A$ and $B$ of a non-abelian group $G$ are said to be
	\begin{enumerate}
		\item{left-homometric if $L_A=L_B$ as multisets.}
		\item{right-homometric if $R_A=R_B$ as multisets.}
	\end{enumerate}
\end{definition}

Obviously, if $B$ is obtained from $A$ by a right-translation, i.e. $B=Ag$ with $g \in G$, then $A$ and $B$ are trivially left-homometric, and similarly if $B=gA$ with $g \in G$ then $A$ and $B$ are trivially right-homometric. As before, left- or right-homometry is understood from now on to always refer to the non-trivial cases (note the map $x \to x^{-1}$ does not give homometric sets in the non-abelian case). Obviously, all notions of abelian or non-abelian homometry can be extended to subsets of sets on which $G$ has a simply transitive action.

In a previous musically-oriented contribution \cite{GenuysPopoff2017}, we have explored left- and right-homometry in dihedral groups. Such groups frequently appear in transformational music theory, in particular the dihedral group $D_{24}$ with 24 elements: it is known as the $T\text{/}I$ group when acting on the left on the set of the 24 major and minor triads, and as the neo-Riemannian $PLR$ group when acting on the right on the same set. If we consider the two sets of chords $A=\{C_m, C\sharp_M, G\sharp_m, E_m, E\flat_M\}$ and $B=\{C_m, E\flat_M, G\sharp_m, F_M, E_m\}$, then the diagrams of intervals in the $T\text{/}I$ group

\begin{center}
	\scalebox{0.68}{
		\TransformationDiagram{$C_m$, $C\sharp_M$, $G\sharp_m$, $E_m$, $E\flat_M$}{
			0/1/$I_{10}$,
			0/2/$T_8$,
			0/3/$I_0$,
			0/4/$T_4$,
			1/0/$I_{10}$,
			1/2/$I_6$,
			1/3/$T_2$,
			1/4/$I_2$,
			2/0/$T_4$,
			2/1/$I_6$,
			2/3/$I_8$,
			2/4/$T_8$,
			3/0/$I_0$,
			3/1/$T_{10}$,
			3/2/$I_8$,
			3/4/$I_4$,
			4/0/$T_8$,
			4/1/$I_2$,
			4/2/$T_4$,
			4/3/$I_4$
		}
	}
	\scalebox{0.68}{
		\TransformationDiagram{$C_m$, $E\flat_M$, $G\sharp_m$, $F_M$, $E_m$}{
			0/1/$I_{8}$,
			0/2/$T_8$,
			0/3/$T_4$,
			0/4/$I_{10}$,
			1/0/$I_{8}$,
			1/2/$I_4$,
			1/3/$I_0$,
			1/4/$T_2$,
			2/0/$T_4$,
			2/1/$I_4$,
			2/3/$T_8$,
			2/4/$I_6$,
			3/0/$T_8$,
			3/1/$I_{0}$,
			3/2/$T_4$,
			3/4/$I_2$,
			4/0/$I_{10}$,
			4/1/$T_{10}$,
			4/2/$I_6$,
			4/3/$I_2$
		}
	}
\end{center}
show that $A$ and $B$ are left-homometric. If we instead consider intervals in the $PRL$ group, i.e. right intervals, we then get the following diagrams

\begin{center}
	\scalebox{0.68}{
		\TransformationDiagram{$C_m$, $C\sharp_M$, $G\sharp_m$, $E_m$, $E\flat_M$}{
			0/1/$LRL$,
			0/2/$PL$,
			0/3/$LP$,
			0/4/$R$,
			1/0/$LRL$,
			1/2/$PRL$,
			1/3/$PRP$,
			1/4/$RLRL$,
			2/0/$LP$,
			2/1/$PRL$,
			2/3/$PL$,
			2/4/$RLP$,
			3/0/$PL$,
			3/1/$PRP$,
			3/2/$LP$,
			3/4/$RPL$,
			4/0/$R$,
			4/1/$RPRL$,
			4/2/$RLP$,
			4/3/$RPL$
		}
	}
	\scalebox{0.68}{
		\TransformationDiagram{$C_m$, $E\flat_M$, $G\sharp_m$, $F_M$, $E_m$}{
			0/1/$R$,
			0/2/$PL$,
			0/3/$PRL$,
			0/4/$LP$,
			1/0/$R$,
			1/2/$RLP$,
			1/3/$RLRL$,
			1/4/$RPL$,
			2/0/$LP$,
			2/1/$RLP$,
			2/3/$PRP$,
			2/4/$PL$,
			3/0/$PRL$,
			3/1/$RPRL$,
			3/2/$PRP$,
			3/4/$LRL$,
			4/0/$PL$,
			4/1/$RPL$,
			4/2/$LP$,
			4/3/$LRL$
			}
	}
\end{center}
which show that $A$ and $B$ are right-homometric. The two sets $A$ and $B$ therefore constitute an example of \textit{simultaneous left- and right-homometry}.

In his seminal book \cite{Lewin1987}, David Lewin considered the positive affine group $\Aff=\{(u,d):u\in\mathbb R,\ d>0\}$ with group operation $(u,d)(v,e)=(u+dv,de)$ which he called the \textit{time-span group}, along with its simply transitive right action on \textit{time-spans}, i.e. pairs $(t,\delta) \in \Aff$ representing the data of a time point on $\mathbb{R}$ together with a positive duration $\delta$. Lewin considers the difference (interval) between two time-spans $(t_1,\delta_1)$ and $(t_2,\delta_2)$ given by $\left(\frac{t_2-t_1}{\delta_1},\frac {\delta_2}{\delta_1}\right)$, that is the difference for the right action of $\Aff$, which he favored for musicological reasons over the left action. It is nevertheless possible to define an interval between two time-spans $(t_1,\delta_1)$ and $(t_2,\delta_2)$ for the left action, which is then equal to $\left(t_2-\frac{\delta_2}{\delta_1}t_1,\frac{\delta_2}{\delta_1}\right)$.

\begin{figure}[t!]
	\centering
	\includegraphics{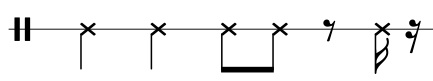}
	\includegraphics{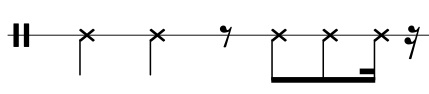}
	\caption{The subsets $\{(0,1),(1,1),(2,1/2),(5/2,1/2),(7/2,1/4)\}$ and $\{(0,1),(1,1),(5/2,1/2),(3,1/2),(7/2,1/4)\}$ of $\Aff$, represented above as rhythms in musical notation, are right-homometric for the action of the positive affine line group.}
	\label{fig:righthomometric_rhythms}
\end{figure}

In a previous contribution on homometric sets in mathematical musicology \cite{Mandereau1,Mandereau2}, Mandereau et al. put forward the two rhythms of Figure \ref{fig:righthomometric_rhythms} corresponding to the two subsets $A=\{(0,1),(1,1),(2,1/2),(5/2,1/2),(7/2,1/4)\}$ and $B=\{(0,1),(1,1),(5/2,1/2),(3,1/2),(7/2,1/4)\}$ of $\Aff$. It can easily be checked that $A$ and $B$ are right-homometric, but not left-homometric, for the simply transitive action of $\Aff$. This leads to the natural question: \textit{are there simultaneous left- and right-homometric subsets of the positive affine line group $\Aff$} ?

The answer to this question is deceptively simple. Let $T=(1,1)$ and $S=(0,1/2)$ in $\Aff$ \footnote{The group generated by $T$ and $S$ is actually isomorphic to the Baumslag-Solitar $\text{BS}(1,2)$ group, with $ST^2=TS$}. For any two elements $TS^m$ and $TS^n$ of $\Aff$, the left difference is $TS^m(TS^n)^{-1} = TS^{m-n}T^{-1}$ while the right difference is $(TS^m)^{-1}TS^n = S^{n-m}$. Left or right differences of these types generate in both cases an infinite cyclic group, and we thus can use the above result from the turnpike problem exposed above: the two subsets $A=\{T,TS,TS^2,TS^6,TS^8,TS^{11}\}$ and $B=\{T,TS,TS^6,TS^7,TS^9,TS^{11}\}$ are simultaneously left- and right-homometric.

The main contribution of this paper is the following answer to the refined question above.

\begin{theorem}\label{thm:main}
	Let $A,B$ be nonempty finite subsets of $\Aff$, such that $L_A=L_B$ and $R_A=R_B$. If the subgroup generated by $R_A$ is nonabelian, then $A=B$. The same conclusion holds if nonabelianity is assumed for $L_A$.
\end{theorem}

The left and right subgroup are actually conjugate so proving the first assertion suffices. The proof proceeds in two steps. A first action of $\Aff$ on a vector space of polynomials shows that the dilation profiles of the subsets $A$ and $B$ are identical. A second action of $\Aff$ on a different vector space of polynomials recovers every individual group elements. Before these two steps are exposed, we prove the following useful lemmas.

\begin{lemma}\label{lem:moments}
	If $d_1,\ldots,d_r$ are distinct positive real numbers and
	$\sum_i c_i d_i^p=0$ for every integer $p\geq1$, then all $c_i$ vanish.
	Also, two finitely supported measures on the real line with the same
	moments of every nonnegative integer order are equal.
\end{lemma}
\begin{proof}
	For the first assertion, use $p=1,\ldots,r$. The matrix $(d_i^p)_{p,i}$
	has nonzero determinant
	$\bigl(\prod_i d_i\bigr)\prod_{i<j}(d_j-d_i)$.
	For the second, subtract the measures and apply the ordinary Vandermonde
	matrix, with powers $0,\ldots,r-1$, to their finite union of supports.
\end{proof}

\begin{lemma}\label{lem:coupling}
	Let $F,G$ be upper triangular real matrices in the ordered basis
	$e_0,\ldots,e_n$, with $n\geq2$. Suppose
	\[
	FG=GF,\qquad F_{00}=F_{nn}=\lambda,\qquad
	F_{jj}\ne\lambda\quad(0<j<n).
	\]
	If $F$ has no $\lambda$-eigenvector with nonzero $e_n$ coefficient, then
	$G_{00}=G_{nn}$.
\end{lemma}
\begin{proof}
	Let $U=\operatorname{span}(e_0,\ldots,e_{n-1})$ and $U'=\operatorname{span}(e_1,\ldots,e_{n-1})$.
	We first show that it is possible to find a vector $u\in U$ such that $(F-\lambda I)(u+e_n)=c e_0$. 
	In the basis $e_0,\ldots,e_n$ the operator $F-\lambda I$ has the matrix form
	\[
F-\lambda I =
\begin{pmatrix}
0 & B & d \\
0 & C & D^{\mathsf T} \\
0 & 0 & 0
\end{pmatrix},
\]
where $C$ is upper triangular, $B$ and $D$ are vectors in $U'$, and $d$ is a scalar.	 Under the assumption of the lemma, $C$ is invertible. Take $x  = u + e_n$ with $u \in U$. We have $u = \alpha e_0 + u'$ with $u' \in U'$, so that
	\[
(F-\lambda I)
\begin{pmatrix}
\alpha\\
u'\\
1
\end{pmatrix}
=
\begin{pmatrix}
Bu'+d\\
Cu'+D^{\mathsf T}\\
0
\end{pmatrix}.
\]
Since $C$ is invertible, it suffices to take $u' = -C^{-1}D^{\mathsf T}$ to obtain
	\[
(F-\lambda I)
\begin{pmatrix}
\alpha\\
-C^{-1}D^{\mathsf T}\\
1
\end{pmatrix}
=
\begin{pmatrix}
-BC^{-1}D^{\mathsf T}+d\\
0\\
0
\end{pmatrix},
\]
i.e. we have $(F-\lambda I)x  = c e_0$ for some $c$. In this case, $c\ne0$, since otherwise $x = u+e_n$ would be a forbidden eigenvector.
	Set $r_0=G_{00}$ and $r_n=G_{nn}$. Triangularity gives
	$w=Gx-r_nx\in U$, and we have
	\begin{align*}
		(F-\lambda I)w &= (F-\lambda I)(Gx-r_nx)\\
					   &= FGx - r_nFx-\lambda Gx + \lambda r_n x\\
					   &= GFx -\lambda Gx - r_n(F - \lambda I)x\\
					   &= G(F - \lambda I)x - r_n c e_0\\
					   &= G c e_0 - r_n c e_0\\
					   &= c(r_0-r_n) e_0\\
	\end{align*}
	Since $F-\lambda I$ induces an invertible operator on $U'$, if the components of the image of $w$ by $(F-\lambda I)$ on $e_1,\ldots,e_{n-1}$ are all null, then the components of $w$ on $e_1,\ldots,e_{n-1}$ are null as well, and $w = k e_0$ for some $k$. Thus $(F-\lambda I)w = 0$, hence $c(r_0-r_n)=0$, proving the claim.
\end{proof}

\section{Reconstructing the dilations}

The first part of the proof consists in reconstructing the dilation profiles of each set. For a given integer $p\geq1$, we consider the vector space $V_p$ of real polynomials of degree at most 2p. We consider the representation $\rho$ of $\Aff$ on $V_p$ which, for a given element $(u,d) \in \Aff$, is such that
\[ (\rho_{(u,d)}P)(x) = d^p P((x-u)/d).\]
It can easily be checked that $\rho$ is indeed a representation, and in the basis $(1,x,\ldots,x^{2p})$ the corresponding matrix $M_{(u,d)}$ is upper triangular with diagonal elements $d^{p-i}$ with $0 \leq i \leq 2p$.

Let $A$ and $B$ be non-empty finite subsets of $\Aff$ such that $L_A=L_B$ and $R_A=R_B$. We form the matrices
\[M_A = \sum_{(u,d) \in A} M_{(u,d)}, \qquad M_B = \sum_{(u,d) \in B} M_{(u,d)},\]
as well as the matrices
\[M_A^\sharp = \sum_{(u,d) \in A} M^{-1}_{(u,d)}, \qquad M_B^\sharp = \sum_{(u,d) \in B} M^{-1}_{(u,d)}.\]
We also form the matrix $F=M_A^\sharp M_A$ and since $A$ and $B$ are left- and right-homometric, we have $F=M_A^\sharp M_A=M_B^\sharp M_B$. 

For each dilation value $d$ occuring in the elements $(u,d)$ of $A$, let $A_d = \{u, (u,d) \in A\}$, and $n_d = |A_d|$. For a given $z \in \mathbb{R}$, we define \[N_A(z) = \sum_d n_d d^z,\] and \[H_A(z) = N_A(z)N_A(-z),\]which corresponds to the sum of exponents of dilation ratios. By homometry we have $H_A(z)=H_B(z)$, which we'll notate as $H$. The diagonal elements of $F=M_A^\sharp M_A=M_B^\sharp M_B$ can be checked to be $H(p-i)$ with $0 \leq i \leq 2p$. In particular we notate $\lambda = H(p) = H(-p)$ and $\lambda > H(p-i)$ for $0 < i < 2p$.

We now prove the following proposition.

\begin{proposition}
	The operator $F$ has no $\lambda$-eigenpolynomial of degree $2p$.
\end{proposition}
\begin{proof}
	Assume there exists a $\lambda$-eigenpolynomial $P$ of degree $2p$ of $F$. The coefficient on the term $x^{2p}$ can be assumed to be positive. The polynomial $P$ thus tends to positive infinity for $x \to +\infty$ and $x \to -\infty$: it therefore has a finite set $S \subset \mathbb{R}$ of minimizers.
	
	Since $P$ is a $\lambda$-eigenpolynomial of $F$, we have
	\[P(x) = \sum_{(u,d) \in R_A} \dfrac{d^p}{\lambda}P((x-u)/d).\]
	For all $x \in S$, $P(x)$ is minimum and so must be $\sum_{(u,d) \in R_A} \dfrac{d^p}{\lambda}P((x-u)/d)$. Since $\sum_{(u,d) \in R_A} \dfrac{d^p}{\lambda} = 1$ and these weights are positive, this means that for all $x \in S$, $P((x-u)/d)$ must be minimum as well. In other words, for each $(u,d) \in R_A$, the map $x \to (x-u)/d$ from $S$ to $S$ is injective. Since each map is also strictly increasing on the ordered set $S$, this means that each is the identity on $S$.
	
	If $S$ has more than one element, then all maps would be the identity on $\mathbb{R}$, and therefore $R_A$ would be trivially abelian, a contradiction.
	If $S$ has only one element $c$, then all maps fix $c$, which imposes $u=c(1-d)$. But then all maps would commute, and the generated subgroup would be abelian, a contradiction.
\end{proof}

Let $G$ be the matrix $G=M_A^{-1}M_B$. We have $FG = M_A^\sharp M_A M_A^{-1}M_B = M_A^\sharp M_B$, and
\[GF = M_A^{-1}M_B M_A^\sharp M_A = M_A^{-1}M_B M_B^\sharp M_B = M_A^{-1}M_A M_A^\sharp M_B = M_A^\sharp M_B
\]
so $F$ and $G$ commute and applying Lemma \ref{lem:coupling} determines that $G_{00}=G_{nn}$, i.e.
\[\frac{N_B(p)}{N_A(p)} = \frac{N_B(-p)}{N_A(-p)}.\]
Since \[\frac{N_B(p)N_B(-p)}{N_A(p)N_A(-p)}=1,\] and since $\frac{N_B(p)}{N_A(p)}$ is positive, we have that $N_A(p)=N_B(p)$ for all $p$. Applying Lemma \ref{lem:moments} we conclude that $A$ and $B$ have the same dilations $d$ and the same counts $n_d$. We write their common function as $N(z)$, so $H(z) = N(z)N(-z)$.

\section{Reconstructing the entire set}

For a given integer $q\geq1$, we now consider the vector space $\tilde{V}_q$ of real polynomials of degree at most q. We consider the representation $\tilde{\rho}$ of $\Aff$ on $\tilde{V}_q$ which, for a given element $(u,d) \in \Aff$, is such that
\[ (\tilde{\rho}_{(u,d)}P)(x) = P((x-u)/d).\]
In the basis $(1,x,\ldots,x^{q})$ the corresponding matrix $\widetilde{M}_{(u,d)}$ is upper triangular with diagonal elements $d^{-i}$ with $0 \leq i \leq q$.

As before, we form the matrices
\[\widetilde{M}_A = \sum_{(u,d) \in A} \widetilde{M}_{(u,d)}, \qquad \widetilde{M}_B = \sum_{(u,d) \in B} \widetilde{M}_{(u,d)},\]
as well as the matrices $\widetilde{F}=\widetilde{M}_A^\sharp \widetilde{M}_A$ and $\widetilde{G} = \widetilde{M}^{-1}_A \widetilde{M}_B$. It is easily checked that $\widetilde{F}\widetilde{G} = \widetilde{G}\widetilde{F}$.

The diagonal elements of $\widetilde{F}$ are $H(i)$ for $0 \geq i \geq q$, which are all distincts and non-zero. Therefore $\widetilde{F}$ is diagonalizable, and $\widetilde{G}$ preserves each of its 1-dimensional eigenspace, since it commutes with $\widetilde{F}$. The diagonal elements of $\widetilde{G}$ are $N_B(i)/N_A(-i)$ for $0 \leq i \leq q$; by the result from the previous section, all diagonal elements are therefore equal to 1. Thus, in the eigenbasis of $\widetilde{F}$, $\widetilde{G}$ is the identity matrix, and we conclude that $\widetilde{M}_A = \widetilde{M}_B$ for all $q \geq 1$. This means that for every polynomial $P$, we have \[\sum_{(u,d) \in A} P((x-u)/d) = \sum_{(u',d') \in B} P((x-u')/d').\]

It is always possible to find $x_0 \in \mathbb{R}$ such that there is a bijection between $A \cup B$ and $\{(x_0-u)/d, (u,d) \in A \cup B\}$. By Lagrange interpolation, we can build a family of polynomials $P_h$ indexed by $h \in A \cup B$ such that $P_h((x_0-u)/d)=1$ if $(u,d)=h$, and $P_h((x_0-u)/d)=0$ if $(u,d) \neq h$. Using $P_h$ in the above equality, we conclude that $(u,d) \in A$ if and only if $(u,d) \in B$, thus proving that $A=B$, which proves Theorem \ref{thm:main}.

\section*{Acknowledgement}
The results presented in this paper originated from extensive guided sessions with GPT-5.1 Sol, and were adapted from it.

\bibliography{affine_nonabelian_homometry}
\bibliographystyle{alpha}
\addcontentsline{toc}{section}{References}

\end{document}